\documentclass[centertags,12pt]{amsart}
\usepackage{latexsym}
\usepackage{amssymb}

\numberwithin{equation}{section}
\makeatletter
\renewcommand{\subsection}{\@startsection
{subsection}{2}{0mm}{\baselineskip}{-0.25cm}
{\normalfont\normalsize\bf}}
\makeatother

\newtheorem*{theorem*}{Theorem}
\newtheorem{theorem}{Theorem}[section]

\newtheorem{corollary}[theorem]{Corollary}
\newtheorem{lemma}[theorem]{Lemma}

\theoremstyle{remark}
\newtheorem{remark}[theorem]{Remark}
\newtheorem{claim}[theorem]{Claim}

\def\A{\mathbf A}

\def\P{\mathbf P}
\def\N{\mathbf N}
\def\Z{\mathbf Z}

\def\cD{\mathcal D}

\def\cX{\mathcal X}

\def\fq{\mathbf F_{q^2}}
\def\bfq{\bar{\mathbf F}_{q^2}}

\def\supp{{\rm Supp}}
\def\div{{\rm div}}
\def\dim{{\rm dim}}
\def\deg{{\rm deg}}

\def\fro{{\mathbf \Phi}}

\begin{document}

\author[M.~Abd\'on]{Miriam Abd\'on} 
\author[F.~Torres]{Fernando Torres}\thanks{{\em 2000 Math. Subj. Class.}: Primary 11G20, Secondary 14G05, 14G10}
\thanks{The authors were partially supported respectively by FAPERJ and PRONEX-Cnpq (Brazil), and by Cnpq-Brazil (Proc. 300681/97-6)}

\title{On a $\mathbf{F_{q^2}}$-maximal curve of genus $\mathbf{q(q-3)/6}$}

\address{Dep. Matem\'atica, PUC-Rio\\
Marqu\^es de S. Vicente 225\\
22453-900, Rio de Janeiro, RJ, Brazil}
\email{miriam@mat.puc-rio.br}

\address{IMECC-UNICAMP, Cx. P. 6065, Campinas, 13083-970-SP, Brazil}
\email{ftorres@ime.unicamp.br}
\address{Current Address: Dpto. Algebra, Geometr\'{\i }a y Topolog\'{\i }a; Fac. Ciencias-Universidad de Valladolid, c/Prado de la Magdalena s/n, Valladolid (Spain)}
\email{ftorres@agt.uva.es}

   \begin{abstract} We show that a $\fq$-maximal curve of genus
$q(q-3)/6$ in characteristic three is either a non-reflexive space curve of degree $q+1$, or it is uniquely determined up to $\fq$-isomorphism by a plane model of Artin-Schreier type.
   \end{abstract}

\maketitle

  \section{Introduction}\label{s1}

Let $\cX$ be a projective, geometrically irreducible, non-singular
algebraic curve of genus $g$ defined over the finite field $\fq$ of order
$q^2$. The curve $\cX$ is called {\em $\fq$-maximal} if it attains the
Hasse-Weil upper bound on the number of $\fq$-rational points; i.e., if
one has
   $$
\#\cX(\fq)=q^2+1+2qg\, .
   $$ 
Maximal curves are known to be very useful in Coding Theory \cite{goppa}
and they have been intensively studied by several authors: see e.g. 
\cite{sti-x}, \cite{vdG-vdV}, \cite{ft1}, \cite{fgt1}, \cite{ft2}, \cite{chkt}, \cite{gt1}, \cite{g-sti-x}, \cite{kt1}, \cite{kt2}. The subject of this paper is related to the following basic questions:
   \begin{itemize}
\item For a given power $q$ of a prime, which is the spectrum of the
genera $g$ of $\fq$-maximal curves?
\item For each $g$ in the previous item, how many non-isomorphic
$\fq$-maximal curves of genus $g$ do exist?
\item Write down an explicit $\fq$-plane model for each of the curves in
the previous item.
   \end{itemize} 
Ihara \cite{ihara} observed that $g$ cannot be arbitrarily large compared with $q^2$. More precisely, 
   $$
g\leq g_1=g_1(q^2):=q(q-1)/2\, .
   $$
R\"uck and Stichtenoth \cite{rueck-sti} showed that (up to
$\fq$-isomorphism) there is just one $\fq$-maximal curve of genus $g_1$,
namely the Hermitian curve of equation
  \begin{equation}\label{eq1.1}
Y^qZ+YZ^q=X^{q+1}\, .
  \end{equation}
On the other hand, if $g<g_1$, then 
   $$
g\leq g_2=g_2(q^2):=\lfloor (q-1)^2/4\rfloor
   $$ 
(see \cite{sti-x}, \cite{ft1}) and, up to $\fq$-isomorphism, there is just
one $\fq$-maximal curve of genus $g_2$ which is obtained as the quotient of the Hermitian curve by a certain involution (see \cite[Thm. 3.1]{fgt1}, \cite{at1}, \cite[Thm. 3.1]{kt2}). Now if $g<g_2$, then (see \cite{kt2})
   $$
g\leq g_3=g_3(q^2):=\lfloor(q^2-q+4)/6\rfloor\, ,
   $$ 
and this bound is sharp as examples in \cite{g-sti-x}, \cite{kt1}, and
\cite{ckt2} show. These examples arise as quotient curves of the Hermitian curve by a certain automorphism of order three; however it is not known whether or not such curves are $\fq$-unique. In view of the results stated above and taking into consideration the examples in \cite{ckt1}, \cite{ckt2} and \cite{g-sti-x}, it is reasonable to expect that only few (non-isomorphic) $\fq$-maximal curves do exist having genus $g$ close to the upper bound $g_1$ provided that $q$ is fixed. As a matter of fact, in the range
  $$
\lfloor (q-1)(q-2)/6\rfloor\leq g <g_3\, ,
  $$
the following statements hold:
   \begin{enumerate}
\item[\rm(I)] If $q\equiv 2\pmod{3}$, there exists an $\fq$-maximal curve
of genus $g=g_3-1$; see \cite[Thm. 6.2]{ckt1} and \cite[Thm. 5.1]{g-sti-x}. 
Such a curve is also the quotient of the Hermitian curve by a certain
automorphism of order three and it is also not known whether this curve is
unique or not;
\item[\rm(II)] If $q\equiv 2\pmod{3}$ and $q\geq 11$, there is just one
$\fq$-maximal curve (up to $\fq$-isomorphism) of genus $(q-1)(q-2)/6$,
namely the non-singular model of the affine plane curve
$y^q+y=x^{(q+1)/3}$, see \cite[Thm. 4.5]{kt2};
\item[\rm(III)] If $q\equiv 1\pmod{3}$ with $q\geq 13$, there is no
$\fq$-maximal curve of genus $(q-1)(q-2)/6$, loc. cit.;
\item[\rm(IV)] If $q=3^t$, $t\geq 1$, there exists an $\fq$-maximal curve
of genus $g=q(q-3)/6$, namely the non-singular model (over $\fq$) of
the affine plane curve 
    \begin{equation}\label{eq1.2}
\sum_{i=1}^{t}y^{q/3^i}=ax^{q+1}\, ,\qquad\text{with $a\in\fq$ such that $a^{q-1}=-1$}\, .
    \end{equation}
We observe that different values of $a$ in (\ref{eq1.2}) define just one $\fq$-maximal curve (up to $\fq$-isomorphism); see Remark \ref{rem4.1}.
    \end{enumerate}
In this paper, motivated by the questions and the results above, we  investigate the uniqueness (up to $\fq$-isomorphism) of the $\fq$-maximal curve in statement (IV). Our main result is Theorem \ref{thm4.1}, where we show that such a curve is either a non-reflexive space curve of degree $q+1$, or it is uniquely determined (up to $\fq$-isomorphism) by the Artin-Schreier relation (\ref{eq1.2}). Unfortunately, we do not know whether or not such a non-reflexive maximal curve can exist (see   Remarks \ref{rem3.1}, \ref{rem3.2}, \ref{rem3.3} and \ref{rem3.4}). 

The key point in the proof of Theorem \ref{thm4.1} is the fact that an  $\fq$-maximal curve of genus $q(q-3)/6$ can be embedded either in $\P^3(\bfq)$ or in $\P^4(\bfq)$ as a curve of degree $q+1$. This fact follows from the ``natural embedding theorem" \cite[Thm. 2.5]{kt1} combined with Castelnuovo's genus formula (\ref{eq2.2}). In the former case, the curve is non-reflexive; see Remark \ref{rem3.3}. In the latter case, the curve is {\em extremal} in the sense that its genus coincide with Castelnuovo's number $c(q+1,4)$; see formula (\ref{eq2.2}). Several relevant properties of extremal curves are known to be true in positive characteristic thanks to Rathmann's paper \cite{rathmann}. For our purpose, the relevant result on extremal curves is a remark due to Accola in \cite{accola}; see Lemma \ref{lemma4.1} here. It follows then that this extremal $\fq$-maximal curve admits a $\fq$-rational point whose Weierstrass semigroup $H$ is generated by $q/3$ and $q+1$ (see Lemma \ref{lemma4.2}); consequently the curve admits a plane model over $\fq$ given by Eq. (\ref{eq4.1}), the so called ``generalized Weierstrass normal form associated to $H$". Further properties of the linear series associated to the embedding of the curve in $\P^4(\bfq)$ (see Section \ref{s2}) and some tedious computations (which are handled as in \cite{at1}) complete the proof of Theorem \ref{thm4.1}.

We point out that several examples of non-isomorphic $\fq$-maximal curves of genus $g\approx q^2/8$ are known; see \cite[Remark 4.1]{chkt} and \cite{abdon}. We also point out that analogous results to our main theorem  have been noticed in \cite{ag}.

As in previous research (see e.g. \cite{fgt1}, \cite{kt2} and the
references therein), the essential tool used here is the approach of St\"ohr and Voloch \cite{sv} to the Hasse-Weil bound applied to the corresponding linear series in the aforementioned ``natural embedding theorem" \cite[Thm. 2.5]{kt1}.

    \section{Preliminaries}\label{s2}

Throughout the paper we assume $q\geq 9$ since the case $q=3$ is trivial. 
As it is known from \cite{ft1}, any $\fq$-maximal curve $\cX$ is equipped
with its $\fq$-canonical linear series; namely, the complete simple
base-point-free linear series
  $$
\cD=\cD_{\cX}:=|(q+1)P_0|\, ,
  $$
where $P_0$ is an arbitrary $\fq$-rational point of $\cX$. Two basic properties of $\cD$ are that it is very ample \cite[Thm. 2.5]{kt1}, and the following linear equivalence of divisors   
\cite[Cor. 1.2]{fgt1}:
   \begin{equation}\label{eq2.1}
qP+\fro(P)\sim (q+1)P_0\, ,\qquad \forall P\in\cX\, ,  
   \end{equation}  
where $\fro=\fro_{q^2}$ is the Frobenius morphism on $\cX$ relative to
$\fq$. In particular, the definition of $\cD$ is independent of $P_0$. To
deal with the dimension $N$ of $\cD$ we use Castelnuovo's genus bound (for curves in projective spaces) which, for a simple linear series $g^r_d$ on $\cX$, bounds from above the genus $g$ of the curve by means of Castelnuovo's number $c(d,r)$; i.e., one has
   \begin{equation}\label{eq2.2} 
g\leq c(d,r):=\frac{d-1-\epsilon}{2(r-1)}(d-r+\epsilon)\, ,
   \end{equation}
$\epsilon$ being the unique integer with $0\leq \epsilon\leq
r-2$ and $d-1\equiv \epsilon\pmod{r-1}$; see \cite{castelnuovo}, \cite[p.
116]{acgh}, \cite[IV, Thm. 6.4]{hartshorne}, \cite[Cor. 2.8]{rathmann}. 
   \begin{lemma}\label{lemma2.1} For a $\fq$-maximal curve of genus
$g=q(q-3)/6,$  $N\in\{3,4\}.$
   \end{lemma}
   \begin{proof} We have that $N\geq 2$ and that $N=2$ if and
only if $\cX$ is the Hermitian curve whose genus is $g=q(q-1)/2$ (see \cite[Thm. 2.4]{ft2}). Therefore $N\geq 3$. If $N\geq 5$, from (\ref{eq2.2}) and the hypothesis on $g$ we would have $q(q-3)/6\leq (q-2)^2/8$; so $q^2\leq 12$, a contradiction.
   \end{proof}
Next based on the St\"ohr and Voloch Theory \cite{sv}, we summarize some
properties on Weierstrass Point Theory and Frobenius Orders with respect
to the linear series $\cD$. Let 
$\epsilon_0=0<\epsilon_1=1<\ldots<\epsilon_N$ and
$\nu_0=0<\nu_1<\ldots<\nu_{N-1}$ denote respectively the $\cD$-orders and
$\fq$-Frobenius orders of $\cD$. For $P\in\cX$, let
$j_0(P)=0<j_1(P)<\ldots<j_N(P)$ be the $(\cD,P)$-orders of $\cD$, and
$(n_i(P):i=0,1,\ldots)$ the strictly increasing sequence that enumerates
the Weierstrass semigroup $H(P)$ at $P$. We have
    $$
0<n_1(P)<\ldots<n_{N-1}(P)\leq q<q+1\leq n_N(P)\, ,
    $$
and $n_N(P)=q+1$ for $P\in\cX(\fq)$ by (\ref{eq2.1}); furthermore,   
$n_{N-1}(P)=q$ for any $P\in\cX$ (\cite[Prop. 1.9]{fgt1}, \cite[Thm.
2.5]{kt2}). We also have the following facts from \cite[Thm. 1.4, Prop.
1.5]{fgt1}:
    \begin{lemma}\label{lemma2.2} 
    \begin{enumerate}
\item[\rm(1)] $\epsilon_N=\nu_{N-1}=q;$
\item[\rm(2)] $\nu_1=1$ if $N\geq 3;$
\item[\rm(3)] $j_1(P)=1$ for any $P;$
\item[\rm(4)] $j_N(P)=q+1$ if $P\in\cX(\fq)$, otherwise $j_N(P)=q;$
\item[\rm(5)] If $P\in\cX(\fq),$ then the $(\cD,P)$-orders
are $n_N(P)-n_i(P)$, $i=0,1,\ldots,N;$
\item[\rm(6)] If $P\not\in\cX(\fq),$ then the elements
$n_{N-1}(P)-n_i(P),$ $i=0,\ldots,N-1,$ are $(\cD,P)$-orders$.$
   \end{enumerate}
   \end{lemma}

    \section{Case $N=3$}\label{s3}

Let $\cX$ be a $\fq$-maximal curve of genus $g=q(q-3)/6$ and let us keep
the notation in Section \ref{s2}. Then, as we saw in Lemma \ref{lemma2.1},
the dimension $N$ of $\cD=\cD_\cX$ is either 3 or 4. In this section we
point out some consequences of the former possibility.
   \begin{lemma}\label{lemma3.1} If $N=3,$ then $\epsilon_2=3.$
   \end{lemma}
   \begin{proof} Let $S$ be the $\fq$-Frobenius divisor associated to $\cD$
(cf. \cite{sv}). Then $\deg(S)=(\nu_1+\nu_2)(2g-2)+(q^2+3)(q+1)$, where
$\nu_1=1$ and $\nu_2=q$ by Lemma \ref{lemma2.2}. For $P\in\cX(\fq)$ it is
known that (loc. cit.)
  $$
v_P(S)\geq j_1(P)+(j_2(P)-\nu_1)+(j_3(P)-\nu_2)=j_2(P)+1\, .
  $$
Moreover, as $j_2(P)\geq \epsilon_2$, the maximality of $\cX$ implies
$\deg(S)\geq (\epsilon_2+1)((q+1)^2+q(2g-2))\, .$
Now, suppose that $\epsilon_2\geq 4$. Then the above inequality becomes
  $$
(q+1)(q^2-5q-2)\geq (2g-2)(4q-1)\, ,
  $$
which is a contradiction with the hypothesis on $g$. Thus we have shown
that $\epsilon_2\in\{2,3\}$. If $\epsilon_2$ were 2, from \cite[Remark
3.3(1)]{ckt1} we would have $g\geq (q^2-2q+3)/6$, which is again a
contradiction with respect to $g$.
   \end{proof}
   \begin{corollary}\label{corollary3.1} If $N=3,$ then $\dim(2\cD)\geq
9.$ 
   \end{corollary}
   \begin{proof} Since $0,1,3,q$ are $\cD$-orders (Lemma
\ref{lemma2.2}), then it is easy to see that $0,1,2,3,4,6,q,q+1,q+3,2q$
are $2\cD$-orders and the result follows.
   \end{proof}
   \begin{corollary}{\rm (\cite[Lemma 3.7]{ckt1})}\label{cor3.2} If $N=3,$
then there exists a $\fq$-rational point $P$ such that $n_1(P)=q-2.$
   \end{corollary}
We conclude this section with some remarks about the possibility $N=3$.
   \begin{remark}\label{rem3.1} (Related with Weierstrass semigroups) From 
Lemma \ref{lemma2.2} and Corollary \ref{cor3.2} there exists
$P\in\cX(\fq)$ such that $n_1(P)=q-2$, $n_2(P)=q$, $n_3(P)=q+1$; i.e.,
the Weierstrass semigroup $H(P)$ contains the semigroup 
  $$
H :=\langle q-2, q, q+1\rangle =\{(q-2)i: i\in \N\}\cup
  \{qi-2(i-1):i\in \Z^+\}
  $$
whose genus (i.e; $\#(\N\setminus H)$) is equal to $(q^2-q)/6$ (see e.g.
\cite[Lemma 3.4]{ckt1}). How can we complete $H$ in order to get $H(P)$?
We have to choose $q/3$ elements from $\N\setminus H$, and it is easy to
see
that such elements must belong to the set
  $$
\{(i(q-1):i=2,\ldots,q/3-1\}\cup\{q^2-5q/3-1, q^2-5q, q^2-5q/3+1\}\, .
  $$
This set contains $q/3+2$ elements, so we have to exclude two
elements from it. Hence we arrive to the following seven possibilities:
    \begin{enumerate}
\item[\rm(i)] $q-1, 2q-2 \notin H(P)$;
\item[\rm(ii)] $q-1 \notin H(P)$ but $2q-2 \in H(P)$; in this case we
have to eliminate one element from the set $\{q^2-5q/3-1, q^2-5q/3,
q^2-5q/3+1\}$; 
\item[\rm(iii)] $q-1, 2q-2 \in H(P)$; in this case we have to eliminate 
two elements from the set $\{q^2-5q/3-1, q^2-5q/3, q^2-5q/3+1\}$. 
    \end{enumerate}
So far, we do not know how an obstruction (for $N$ being equal to 3) might
arise from some of the possibilities above.
    \end{remark}
   \begin{remark}\label{rem3.2} (Related with the Hermitian curve
(\ref{eq1.1})) Suppose that $\cX$ is $\fq$-covered by the Hermitian curve.
Then the covering cannot be Galois; otherwise by \cite[Prop. 5.6]{ckt2}
the curve would be $\fq$-isomorphic to the non-singular model of
(\ref{eq1.2}) and thus $N=4$. We recall that there is not known any
example of a $\fq$-maximal curve $\fq$-covered by the Hermitian curve by a
non-Galois covering.
   \end{remark}
   \begin{remark}\label{rem3.3} (Reflexivity, Duality and the Tangent Surface) Recall that we can assume our curve $\cX$ as being embedded in
$\P^3(\bfq)$ by \cite[Thm. 2.5]{kt1}. Hefez \cite{hefez1} noticed that
four cases for the generic contact orders for space curves can occur.
Homma \cite{homma2} realized that all the four aforementioned cases occur
and characterized each of them by means of the reflexivity of either the
curve $\cX$, or the tangent surface ${\rm Tan}(\cX)$ associated to it (see
also \cite{homma1}). In our situation ($\epsilon_2=3$ and $\epsilon_3=q$),
the maximal curve $\cX$ is non-reflexive by Hefez-Kleiman Generic Order of Contact Theorem \cite{hefez-klm}. In addition, Homma's result implies that ${\rm Tan}(\cX)$ is also non-reflexive. So far
we do not know how to relate the maximality of $\cX$ to the
non-reflexivity of its tangent surface. We mention that techniques
analogous to those of \cite{sv} that work on certain surfaces in $\P^3$ 
over prime fields is now available thanks to a recent paper by Voloch
\cite{voloch}.
    \end{remark}
    \begin{remark}\label{rem3.4} (Related to Halphen's theorem) Ballico
\cite{ballico} extended Harris \cite{harris} and Rathmann \cite{rathmann}
results concerning space curves contained in surfaces of certain degree.
For a $\fq$-maximal curve $\cX$ of genus $q(q-3)/6$, with $q$ large
enough, Ballico's result implies that $\cX$ is contained in a
surface of degree 3 or 4. On the other hand, suppose that Voloch's
approach \cite{voloch} can be extended to cover the case of surfaces over
arbitrary finite fields. Then a conjunction of Ballico and Voloch's result
would provide further insights of the curves studied here.
    \end{remark}

   \section{Main Result}\label{s4}
Let $\cX$ be a $\fq$-maximal curve of genus $g=q(q-3)/6$ and let us keep the notation in Section \ref{s2}. We will assume that $\cX$ is embedded in $\P^N(\bfq)$ as a curve of degree $q+1$ (cf. \cite[Thm. 2.5]{kt1}), where  $N=\dim(\cD)\in\{3,4\}$ by Lemma \ref{lemma2.1}. 
   \begin{theorem}\label{thm4.1} We have that either
   \begin{enumerate}
\item[\rm(1)] $N=3$ and $\epsilon_2=3;$ i$.$e$.$ the curve is a non-reflexive space curve$;$ or
\item[\rm(2)] $N=4$ and $\cX$ is $\fq$-isomorphic to the non-singular model of the plane curve (\ref{eq1.2})$.$ 
   \end{enumerate}
   \end{theorem}
   \begin{remark}\label{rem4.1} In Case (2), the curve is $\fq$-covered by the Hermitian curve (\ref{eq1.1}). To see this we observe that via the automorphism $(x,y)\mapsto (x,a^qy)$, with $a^{q-1}=-1$, the Hermitian curve admits a plane model of type $y_1^q-y_1=ax^{q+1}$. Then we obtain the curve (\ref{eq1.2}) via the rational function $y:=y_1^3-y_1$. In particular, the non-singular model of (\ref{eq1.2}) is indeed a $\fq$-maximal curve by \cite[Prop. 6]{lachaud}.

Furthermore, the non-singular model of (\ref{eq1.2}) does not depend on the element $a$ used to define the afin equation. In fact if $a_1, a_1\in \fq$ define curves as in (\ref{eq1.2}), then such curves are $\fq$-isomorphic by means of the automorphism $(x,y)\mapsto (\alpha x,y)$ with $\alpha^{q+1}=a_1/a_2$.
   \end{remark}
To give the proof of the theorem we need some auxiliary results. 

Let us assume $N=4$. In this case, the curve $\cX$ becomes extremal in the sense mentioned in the introduction. The following result is implicitly contained in the proof of Castelnuovo's genus bound (\ref{eq2.2}) taking into account the Riemman-Roch theorem; see e.g. \cite[p. 361 and Lemma
3.5]{accola}.
   \begin{lemma}\label{lemma4.1} 
    \begin{enumerate}
\item[\rm(1)] $\dim(2\cD)=11;$
\item[\rm(2)] There exists a base-point-free $2$-dimensional complete
linear series $\cD'$ of degree $2q/3$ such that $\frac{q-6}{3}\cD+\cD'$ is
the canonical linear series of $\cX.$
    \end{enumerate}
    \end{lemma}
Part (1) of this lemma implies Lemmas 4.2 and 4.4 in \cite{kt2}:
   \begin{corollary}\label{cor4.1}\begin{enumerate}
  \item[\rm(1)] If $j_2(P)=2,$ then $j_3(P)=3;$
  \item[\rm(2)] If $P\in\cX(\fq)$ and $j_2(P)>2,$ then $j_2(P)=(q+3)/3,$
$j_3(P)=(2q+3)/3;$
  \item[\rm(3)] If $P\not\in\cX(\fq)$ and $j_2(P)>2,$ then either
$j_2(P)=q/3,$ $j_3(P)=2q/3;$ or $j_2(P)=(q-1)/2),$ $j_3(P)=(q+1)/2.$
    \end{enumerate}
   \end{corollary}
   \begin{lemma}\label{lemma4.11} For $q=9$, $n_1(P)=3$ for any
$P\in\cX(\fq).$
    \end{lemma}
    \begin{proof} 
Let $P$ be a $\fq$-rational and set $n_i:=n_i(P)$. Lemma \ref{lemma2.2}
implies that $0,1,j_2=10-n_2,j_3=10-n_1,10$ are the $\cD$-orders at the
point. Then the set of $2\cD$-orders at $P$ must contain the set
$\{0,1,2,j_2,j_2+1,2j_2,j_3,j_3+1,j_2+j_3,2j_3,10,11,j_2+10,j_3+10,20\}$
and hence, as dim$(2\cD)=11$ by Lemma \ref{lemma4.1}, the result follows. 
    \end{proof}
From now on let us assume $q\geq 27.$
   \begin{lemma}\label{lemma4.2}\begin{enumerate}
\item[\rm(1)] The case (3) in Corollary \ref{cor4.1}
cannot occur$;$
\item[\rm(2)] There exists $P_1\in\cX(\fq)$ such that $j_2(P_1)>2;$ in
this case, $n_1(P_1)=q/3;$
\item[\rm(3)] Let $P_1$ be as in (2) and $x\in\fq(\cX)$ such that
$\div_\infty(x)=\frac{q}{3}P_1.$ Then the morphism
$x:\cX\setminus\{P_1\}\to\A^1(\bfq)$ is unramified$,$ and $x^{-1}(\alpha)\subseteq \cX(\fq)$ for any $\alpha\in\fq;$
\item[\rm(4)] The $\cD$-orders and $\fq$-Frobenius orders of $\cD$ are
respectively $0,1,2,3,q$ and $0,1,2, q.$
    \end{enumerate}
    \end{lemma}
    \begin{proof} For $P\in\cX$, set $j_i=j_i(P)$ and $n_i=n_i(P)$. 

(1) We have that $\{q-n_2,q-n_1\}\subseteq \{1,j_2,j_3\}$ by Lemma
\ref{lemma2.2}(6). Suppose that Case (3) in Corollary. \ref{cor4.1} occurs.

{\em Case $j_2=q/3$, $j_3=2q/3$.} Here $n_1\in\{2q/3,q/3\}$; let
$f\in\bfq(\cX)$ such that $\div(f-f(\fro(P)))=D+e\fro(P)-n_1P$, where
$e\geq 1$ and $P\not\in\supp(D)$. If $n_1=2q/3$, then $3e+2$ is an
$(2\cD,\fro(P))$-order by (\ref{eq2.1}). However, as $\dim(2\cD)=11$, the
sequence of $(2\cD,\fro(P))$-orders is
$0,1,2,q/3,q/3+1,2q/3,2q/3+1,q,q+1,4q/3,5q/3,2q$ and thus $n_1=q/3$. In
this case, arguing as above, $3e+1$ is an $(\cD,\fro(P))$-order which is a
contradiction.

{\em Case $j_2=(q-1)/2$, $j_3=(q+1)/2$.} From Lemma
\ref{lemma4.1}(2) and the hypothesis $q\geq 27$, we have that $2j_2+1=q$
is a Weierstrass gap at $P$ (i.e.; $q\not\in H(P)$), a contradiction with
$n_3=q$.

(2) If we show that there exists $P_1\in\cX$ such that $j_2(P_1)>2$, then
the point $P_1$ will be $\fq$-rational by (1). So, suppose that
$j_2(P)=2$ for any $P\in\cX$. Let $R$ denote the ramification divisor
associated to $\cD$ (cf. \cite{sv}). Then the $\cD$-Weierstrass points
coincide with the
set of $\fq$-rational points and $v_P(R)=1$ for $P\in\cX(\fq)$ (cf. Lemma
\ref{lemma2.2}). Therefore
   $$
\deg(R)=(q+6)(2g-2)+5(q+1)=\#\cX(\fq)=(q+1)^2+q(2g-2)\, ,
   $$
so that $2g-2=(q-1)(q-4)/6$, a contradiction. That $n_1=q/3$ follows
immediately from Corollary. \ref{cor4.1}(2) and Lemma \ref{lemma2.2}(5).
    
(3)-(4) By (\ref{eq2.1}) we can assume that $P_0$ satisfies Item (2). Let  $x, y\in\fq(\cX)$ be such that $\div(x)=\frac{q}{3}P_0$ and $\div_\infty(y)=(q+1)P_0$. Then the sections of $\cD$ are generated
by $1,x,x^2,x^3,y$. Now, if we show that there exists $P\in\cX(\fq)$ such that $j_2(P)=2$ and $j_3(P)=3$, then
(4) follows since $\epsilon_i\leq j_i(P)$ and
$\nu_{i-i}\leq j_i(P)-j_1(P)$ (cf. \cite{sv}). To see that such a point
$P$ does exist, we proceed as in \cite[p. 38]{ft2}. For $P\in\cX\setminus\{P_0\}$, write
$\div(x-x(P))=eP+D-n_1P_0$ with $e\geq 1$ and $P,P_0\not\in\supp(D)$. Then
$e,2e,3e$ are $(\cD,P)$-orders and if $e>1$, $3e=q+1$, a contradiction as
$q\equiv 0\pmod{3}$. Finally, for $\alpha\in\fq$ we show that $x^{-1}(\alpha)\subseteq \cX(\fq)$. Let $a_\alpha:=\#x^{-1}(\alpha)$. Then 
   $$
\#\cX(\fq)=1+\frac{q^3}{3}=1+\sum_{\alpha\in \fq}a_\alpha
   $$
and hence $a_\alpha\geq q/3$ for any $\alpha$.
    \end{proof}
In what follows we will assume that $P_0$ satisfies Lemma \ref{lemma4.2}(2)(3); cf. (\ref{eq2.1}). Let $x,y\in\fq(\cX)$ such that $\div_\infty(x)=\frac{q}{3}P_0$ and $\div_\infty(y)=(q+1)P_0$. Let $v=v_{P_0}$ be the valuation at $P_0$, and $D^i:=D^i_x$ the $i$-th
Hasse differential operator on $\bfq(\cX)$ with respect to $x$ (see e.g.
\cite[\S3]{hefez}). We set $D:=D^1$. 
    \begin{corollary}\label{cor4.2}\quad $v(Dy)=v(x^q)=-q^2/3.$
    \end{corollary}
    \begin{proof} (cf. \cite[p. 47]{at1}) Let $t$ be a local parameter at
$P_0$; then 
     $$
v(Dy)=v(dy/dt)-v(dx/dt)=-q-2-(2g-2)=-q^2/3\, ,
     $$
since by the previous lemma the morphism $x:\cX \to \P^1(\bfq)$ is
totally ramified at $P_0$ and unramified outside $P_0$ so that $\div(dx)=(2g-2)P_0$.
     \end{proof}
In what follows, we will use several times the following basic property of $v$: For $f_1, \ldots, f_m\in\bfq(\cX)$ we have that
   \begin{equation}\label{eq4.0}
v(\sum_{i=1}^mf_i)=\min\{v(f_i):i=1,\ldots,m\}\, ,\qquad\text{provided that $v(f_i)\neq v(f_j)$ for $i\neq j$}\, .
   \end{equation}
{\bf Proof of Theorem \ref{thm4.1}.} We already know that either $N=3$ or $N=4$. If $N=3$ the result follows from Remark \ref{rem3.3}.

Let $N=4$. The rational functions $x$ and $y$ above are related to each other by an equation over $\fq$ of type (see e.g. \cite{kato})
   \begin{equation}\label{eq4.1}
x^{q+1}+\sum_{i=0}^{q/3}A_i(x)y^i=0\, ,
   \end{equation}
where the $A_i(x)$'s are polynomials in $x$ such that $\deg(A_i(x))\leq q-3i$, and $A_{q/3}(x)=A_{q/3}\in\fq^*$. In addition by Lemma \ref{lemma4.2}(4) we have a relation of type (cf. \cite[Prop. 2.1]{sv}):
   \begin{equation}\label{eq4.2}
y^{q^2}-y=(x^{q^2}-x)Dy+(x^{q^2}-x)^2D^2y+(x^{q^2}-x)^3D^3y\, .
   \end{equation}
 \begin{claim}\label{claim4.1} For $j=0,\ldots,t-2$ and $k\not\equiv\pmod{3},$ $k\geq 2$ and $3^jk\leq q/3$ we have that $A_{3^jk}(x)=0$ and $A_{3^j}(x)=A_{3^j}\in \fq^*.$
   \end{claim}
(Recall that $q=3^t$.)
   \begin{proof} We apply induction on $j$. First we show that $A_k(x)=0$ if $k\not\equiv 0\pmod{3}$ and $2\leq k\leq q/3$. To do that, apply $D=D^{3^0}$ to Eq. (\ref{eq4.1}); so
  $$
x^q+F_0+G_0Dy=0\, ,\qquad\text{with}
  $$
  $$
F_0:=\sum_{i=0}^{q/3-1}y^iDA_i(x)\, ,\qquad\text{and}\qquad G_0:=(\sum_{i=1}^{q/3-1}A_i(x)iy^{i-1})Dy\,  .
  $$
Suppose that $A_k(x)\neq 0$ for some $k\not\equiv
0\pmod{3}$ and $2\leq k\leq q/3$. Then $v(G)<v(x^q)$ by Corollary \ref{cor4.2} and thus $v(F)=v(GDy)$. Therefore by property (\ref{eq4.0}) 
there exist integers $2\leq i_0\leq q/3-1$, $i_0\not\equiv
0\pmod{3}$ and $1\leq j_0\leq q/3-1$ such that
$v(A_{i_0}(x)y^{i_0-1})=v(y^{j_0}DA_{j_0}(x))$, a contradiction as $\gcd(q/3,q+1)=1$. The same argument shows that $A_1(x)=A_1\in\fq^*$ and thus the induction for $j=0$ is complete.

Let $j=1$. The case $j=0$ reduces Eq. (\ref{eq4.1}) to
   \begin{equation}\label{eq4.3}
x^{q+1}+A_0(x)+A_1y+\sum_{i=1}^{q/9}A_{3i}(x)y^{3i}=0\, .
 \end{equation}
Apply $D^3$ to this equation. We find that
    $$
D^3A_0(x)+A_1D^3y+F_1+G_1(Dy)^3\, ,
    $$
with 
    $$
F_1:=\sum_{i=0}^{q/3}y^{3i}DA_{3i}(x)\, ,\qquad\text{and}\qquad 
G_1:=\sum_{i=0}^{}A_{3i}(x)iy^{3i-3}\, .
    $$
At this point we notice that $v(D^3y)=-q^2$. Indeed, from Eq. (\ref{eq4.2}) we have that $v(D^2y-(x^{q^2}-x)D^3y)=-q^3/3-q^2$ and so it is enough to 
see that $v(D^2y)>-q^3/3-q^2$. This inequality easily follows by applying $D^2$ to Eq. (\ref{eq4.3}) and by using property (\ref{eq4.0}) once more. Now the case $j=1$ follows as in the previous case.

Finally, the case $j$ implies the case $j+1$ (with $j\geq 1$) is even easier than the previous cases since $D^{3^{j+1}}y=0$ by Lemma \ref{lemma4.2}(4).
  \end{proof} 
Therefore Eq. (\ref{eq4.1}) becomes
   \begin{equation}\label{eq4.4}
x^{q+1}+A_0(x)+\sum_{i=0}^{t-1}A_{3^i}y^{3^i}=0\, ,
   \end{equation}
where $A_0(x)$ is a polynomial in $x$ of degree at most $q$, and where each $A_{3^i}\in \fq^*$. Let $A_0(x)=\sum_{i=0}^qB_ix^i$. By using maps of type $(x,y)\mapsto (x+a,y)$ and $(x,y)\mapsto (x, ax+y)$ we can assume that $B_1=B_q=0$. We can also assume that $B_0=0$ by Lemma \ref{lemma4.2}(3). Moreover, it is easy to see that the property: $D^{3^i}y=0$ for $i=2,\ldots,t-1$ (cf. Lemma \ref{lemma4.2}(4)), reduces $A_0(x)$ to 
    $$
A_0(x)=B_2x^2+\sum_{i=1}^{t-1}B_{3^i}x^{3^i}+ 
\sum_{i=1}^{t-1}B_{2\cdot 3^i}x^{2\cdot 3^i}\, .
    $$
Next from Eq. (\ref{eq4.4}) we see that the rational functions $Dy$, $D^2y$  and $D^3y$ are polynomials in $x$. From Eq. (\ref{eq4.2}) we conclude then that $y^{q^2}-y$ is a polynomial in $x$ which is a multiple of $x^{q^2}-x$ by Lemma \ref{lemma4.2}(3). Therefore 
    \begin{itemize}
\item $B_i=0$ for each $i$, and
\item $A_{3^i}=(\frac{A_3}{A_1^3})^{(3^i-1)/2}A_1^{3^i}$ for $i=1, \ldots,t-1$.
    \end{itemize}
Now it is easy to see that there exists $a\in\fq$ such that
$A_3/A_1^3=a^2$ and such that $a^{q-1}=-1$. Finally, we complete the
proof of Theorem \ref{thm4.1} via the map $(x,y)\mapsto (x, aA_1y)$.
    
   \end{document}